\documentclass[author,12pt]{elsarticle}
\usepackage{tikz}
\usepackage{graphicx}
\usepackage{float}
\usepackage{amsthm,epsfig,epstopdf}
\usepackage{amsmath,amssymb,amsfonts}
\usepackage{bbm}
\usepackage[dvipdfm,pdfborder=000,colorlinks,citecolor=blue,linkcolor=blue]{hyperref}
\newtheorem{theorem}{\bf Theorem}[section]

\newtheorem{definition}[theorem]{\bf Definition}
\newtheorem{lemma}{Lemma}[section]
\newtheorem{remark}{Remark}[section]
\newtheorem{example}{Example}[section]

\allowdisplaybreaks

\begin{document}

\renewcommand{\thefootnote}{\fnsymbol{footnote}} 

\renewcommand{\thepage}{\arabic{page}}
\setcounter{page}{1} \begin{center} {\rm \bf \Large Optimal reinsurance and dividends with transaction costs and taxes under thinning structure}\end{center} \vspace{0.1cm}
\begin{center} { Mi Chen$^a$, Kam Chuen Yuen$^b$, Wenyuan Wang$^{c,}$\footnote{Corresponding author.\\ \indent~ \textit{E-mail address:}
wwywang@xmu.edu.cn (W.Y. Wang).}}
\end{center}
\begin{center}
{\small $^a$ College of Mathematics and Informatics \& FJKLMAA, Fujian Normal University,}\\
{\small Fuzhou 350117, China}\\
{\small $^b$ Department of Statistics and Actuarial Science, The University of Hong Kong,}\\
{\small Pokfulam Road, Hong Kong}\\
{\small $^c$ School of Mathematical Sciences, Xiamen University, Xiamen 361005, China}\\
\end{center}

\baselineskip=18pt

\noindent {\sl Abstract:}
In this paper, we investigate the problem of optimal strategies of dividend and reinsurance under the Cram\'{e}r-Lundberg risk model embedded with the thinning-dependence structure which was firstly introduced by Wang and Yuen (2005), subject to the optimality criteria of maximizing the expected accumulated discounted dividends paid until ruin. To enhance the practical relevance of the optimal dividend and reinsurance problem, non-cheap reinsurance is considered and transaction costs and taxes are imposed on dividends, which converts our optimization problem into a mixed classical-impulse control problem.
For the purpose of better mathematical tractability and neat, explicit solutions of our control problem, instead of the Cram\'er-Lundberg framework we study its approximated diffusion model with two thinly dependent classes of insurance businesses.
Using a method of quasi-variational inequalities, we show that the optimal reinsurance follows a two-dimensional excess-of-loss reinsurance strategy, and, the optimal dividend strategy turns out to be an impulse dividend strategy with an upper and a lower barrier, i.e., every thing above the lower barrier is paid as dividends each time the surplus is above the upper barrier, otherwise no dividends are paid.
Closed-form expression for the value function associated with the optimal dividend and reinsurance strategy is also given.
In addition, some numerical examples are presented to illustrate the optimality results.
\\


 \baselineskip=18pt
 \section{Introduction}
  \setcounter{equation}{0}

\indent

As for listed insurance companies, distribution of dividends is a main approach to share profits with policy holders, while purchase of reinsurance is an effective way to reduce risk exposure. Due to the importance of these two features, risk models with reinsurance and dividend payments have received extensive attention in the actuarial literature in the past few decades.  Optimal dividend problem under the diffusion risk model was first investigated by Jeanblanc-Picqu{\'e} and Shiryaev (1995) using the technique of stochastic control theory. Since then, optimal dividend and/or reinsurance problems were studied for different risk models with various objective functions.
There are some well-known dividend strategies that turned out to be optimal in certain situations. For instance,
H{\o}jgaard and Taksar (1999) showed that the optimal dividend strategy is a threshold strategy if the rate of dividend payout is bounded by some positive constant, while it is a barrier strategy for the case where there is no restriction on the rate of dividend payout. When transaction costs is considered, the optimal dividend strategy is usually an impulse strategy, see for example Paulsen (2007, 2008).  The extensively studied risk models for the optimal dividend problem in the literature include diffusion model, Cram\'{e}r-Lundberg model, jump-diffusion model and L\'{e}vy risk model. For example, Asmussen and Taksar (1997), H{\o}jgaard and Taksar (1999), Asmussen et al. (2000), Paulsen (2003), Gerber and Shiu (2004), L{\o}kka and Zervos (2008), He and Liang (2008), Bai et al. (2010), Chen et al. (2013), Yao et al. (2014, 2016), Peng et al. (2016), Vierk\"{o}tter and Schmidli (2017), Zhu (2017), and Liang and Palmowski (2018) considered the optimal dividend problem in the diffusion model; H{\o}jgaard (2002), Azcue and Muler (2005), Schmidli (2006), Gerber and Shiu (2006), Albrecher and Thonhauser (2008), and Azcue and Muler (2012) studied the optimal dividend strategy under the Cram\'{e}r-Lundberg model.
As for other risk models such as the jump-diffusion model and the L\'{e}vy risk model, recent related research can be found in
Avram et al. (2007, 2015), Kyprianou and Palmowski (2007), Loeffen (2008, 2009), Loeffen and Renaud (2010),
Czarna and Palmowski (2010), Wang and Hu (2012), Hunting and Paulsen (2013), Hernandez and Junca (2015), Zhao et al. (2017), P\'{e}rez et al. (2018), Wang et al. (2018), Wang and Zhou (2018), Wang and Zhang (2019), etc.
For more works we are referred to the survey paper Albrecher and Thonhauser (2009) and the references therein.

Although most of the research in this direction mainly deals with independent risks, much attention has been paid to the optimization problems in relation to dependent risks in recent years.  For the risk model with common shock dependence, Bai et al. (2013) derived the optimal excess-of-loss reinsurance strategies that minimize ruin probability; Yuen et al. (2015) and Liang and Yuen (2016) considered the optimal proportional reinsurance strategy under the criterion of maximizing the expected exponential utility; Zhang and Liang (2017) studied the problem of portfolio optimization for jump-diffusion risky assets with common shock dependence and state dependent risk aversion; and Li et al. (2016) investigated the optimal dividend and reinsurance problem in the approximated diffusion model. In recent years, this kind of optimality study has been extended to the risk model with the thinning-dependence structure proposed by Wang and Yuen (2005) which embraces the common shock risk model. Such a generalization undoubtedly makes the problem of study more complicated and challenging.  For example, under the thinning dependence, Han et al. (2018) used the technique of HJB equation to investigate the optimal proportional reinsurance problem that minimize the probability of drawdown in the Brownian motion case; and  Wei et al. (2018) derived the optimal proportional reinsurance strategy in the compound Poisson case under the criterion of maximizing the adjustment coefficient.

In this paper, the problem of optimal dividends and reinsurance under the thinning-dependence structure is studied. We adopt the expected value premium principle and take into account dividend payments subject to transaction costs and taxes.
In order to make our problem mathematically tractable and to obtain neat and explicit solutions for the optimal dividend and reinsurance strategy and its associated value function, instead of the Cram\'er-Lundberg framework we study its approximated diffusion model with two thinly dependent classes of insurance business.
Under this approximated diffusion setup with thinning dependence, we first show that the optimal reinsurance does not have the form of proportional reinsurance strategy that was studied in Han et al. (2018) and Wei et al. (2018), but follows the excess-of-loss reinsurance strategy.  Since fixed transaction costs of dividends are considered, the optimization problem becomes a mixed classical-impulse stochastic control problem, and hence the methods used in Han et al. (2018) and Wei et al. (2018) can not be applied.  By the method of quasi-variational inequalities (QVI), closed-form expressions for the value function and the corresponding optimal excess-of-loss reinsurance and impulse dividend strategy are derived.

Although there are a lot of existing works on the topic of optimal dividend and reinsurance, the literature that takes into account the dependence structure is still fairly scarce. As far as the authors know, in addition to Li et al. (2016), this paper represents the only other attempt in discussing the optimal dividend and reinsurance problems under risk models involving dependence structure.
Compared with the optimal dividend and reinsurance problem without dependent risk, the optimal reinsurance strategy in this paper is a two-dimensional excess-of-loss reinsurance strategy, and the two coordinate-reinsurance strategies are related with each other complicatedly. In order to determine the optimal two-dimensional reinsurance strategy explicitly, we need to define three auxiliary functions and analyse two zeros associated with these auxiliary functions. The optimal dividend and reinsurance control problem is then solved corresponding to two opposite scenarios of the relation of the two zeros.
Compared with Li et al. (2016) that considered the optimal dividend and reinsurance problem with dependent risk too, our paper is quite different in that: Firstly, the common shock dependence structure discussed in Li et al. (2016) is a special case of the thinning-dependence structure of the present paper; Secondly, Li et al. (2016) studied a classical control problem with no transaction costs and taxes, while transaction costs and taxes are considered in this paper which converts our problem into an impulse control problem, and hence different approach as the QVI method is employed; Thirdly, we find the optimal reinsurance strategy dominating all admissible reinsurance strategies to be a particular two-dimensional excess-of-loss reinsurance strategy, while Li et al. (2016) characterized the optimal reinsurance strategy only among the sub-class of excess-of-loss reinsurance strategies.

The rest of this paper is organized as follows. In Section 2, the model and mathematical formulation of the problem are introduced. In Section 3, we show that the excess-of-loss reinsurance strategy is the optimal reinsurance form for our optimization problem. In Section 4, the QVI and verification theorem are presented. Section 5 is devoted to the derivation of the solution to the QVI. The value function and the optimal strategy are given in Section 6. Finally, some numerical examples are provided in Section 7.

\section{The Model}\label{sec:prob-form}
  \setcounter{equation}{0}

\indent We assume that all stochastic quantities are defined on a large enough complete probability space $(\Omega, \mathcal{F},\mathcal{F}_t\,, \mathbf{P})$, where the filtration $\mathcal{F}_t$ represents the information available at time $t$, and any decision made is based on this information.

\indent
The thinning-dependence structure considered in this paper was first introduced by Wang and Yuen (2005). Suppose that an insurance company has a portfolio of $n~(n\geq 2)$ dependent classes of insurance business, and the stochastic sources that may cause a claim in at least one of the classes are classified into $m$ groups. It is assumed that each event occurred in the $k$th group may cause a claim in the $l$th class with probability $p_{kl}$ for $k=1,2,\ldots,m$ and $l=1,2,\ldots,n$, and that for each $l$, there exists at least some $k$ such that $p_{kl}>0$. For the $k$th group, let $N^k(t)$ be the number of events occurred up to time $t$, and $N^k_l(t)$ be the number of claims of the $l$th class up to time $t$ generated from the events in group $k$. For the $l$th class, let $X^{(l)}_i$ ($i=1,2,\ldots$) be the claim size random variables following a common distribution $F_l$ (corresponding to a random variable $X_l$), and denote by $\mu_l$ and $\sigma_l^2$ the mean and the variance of the distribution $F_l$, respectively. Then the aggregate claims process of the company is given by $$ S(t)=\sum_{l=1}^n S_l(t)=\sum_{l=1}^n \sum_{i=1}^{N_l(t)} X^{(l)}_i,$$ where $\{X^{(l)}_i; i=1,2,\ldots\}$ is a sequence of i.i.d. non-negative random variables for each $l$, and $N_l(t)=N^1_l(t)+N^2_l(t)+\cdots+N^m_l(t)$ is the claim-number process of the $l$th class. As usual, we assume that the processes $N^1(t),\ldots,N^m(t)$ are independent Poisson processes with parameters $\lambda_1,\ldots,\lambda_m$, respectively. Furthermore, for $k\neq j$, the two vectors of claim-number processes, $(N^k(t),N^k_1(t),\ldots,N^k_n(t))$ and $(N^j(t),N^j_1(t),\ldots,N^j_n(t))$ are independent; and for each $k$, $N^k_1(t),\ldots,N^k_n(t)$ are conditionally independent given $N^k(t)$. Also, we assume that the $n$ sequences $\{X^{(1)}_i; i=1,2,\ldots\},\ldots,\{X^{(n)}_i; i=1,2,\ldots\}$ are mutually independent and are independent of all the claim-number processes.

The reserve process of the insurer without reinsurance is given by
$$U_t=x+ct-S(t),$$
where $x\geq 0$ is the initial reserve, and $c>0$ is the premium rate.
In order to manage the underlying insurance risk properly, the insurer would like to buy reinsurance to alleviate the impact of large losses. Suppose that the reinsurance strategy for the $l$th class is $q_l$ (not time-varying) with $0\leq q_l(x)\leq x$ for $x\geq 0$ and $l=1,2,\ldots,n$,
and the reinsurance premium rate is denoted by $\delta(\textbf{q})$ with $\textbf{q}=(q_{1},q_{2},\ldots,q_{n})$. Then the reserve process after reinsurance can be written as
$$U^\textbf{q}_t=x+[c-\delta(\textbf{q})]t-S^\textbf{q}(t),$$
where $$ S^\textbf{q}(t)=\sum_{l=1}^n \sum_{i=1}^{N_l(t)} q_{l} (X^{(l)}_i).$$

Similar to Wang and Yuen (2005), we know that $S^\textbf{q}(t)$ follows a compound Poisson process with
\begin{eqnarray*}
 E[S^\textbf{q}(t)]&=& \sum_{l=1}^n E[q_{l} (X^{(l)})]  \sum_{k=1}^m \lambda_k p_{kl}t, \\
 Var[S^\textbf{q}(t)]&=&\sum_{l=1}^n E[q_{l} (X^{(l)})]^2 \sum_{k=1}^m \lambda_k p_{kl}t\\
&&~~+\sum_{l=1}^n \sum_{j\neq l}^n E[q_{j} (X^{(j)})]E[q_{l} (X^{(l)})] \sum_{k=1}^m \lambda_k p_{kj}p_{kl}t.
\end{eqnarray*}
Then $U^\textbf{q}_t$ can be approximated by a pure diffusion $X_t^\textbf{q}$, which is given by
\begin{eqnarray*}
X_t^\textbf{q}=x+ [c-\delta(\textbf{q})-a(\textbf{q})] t+ b(\textbf{q}) W_t,
\end{eqnarray*}
where $\{W_t, t\geq 0\}$ is a standard Brownian motion and
\begin{eqnarray*}
&& a(\textbf{q})= \sum_{l=1}^n E[q_{l} (X^{(l)})]  \sum_{k=1}^m \lambda_k p_{kl}, \\
&& b^2(\textbf{q})=\sum_{l=1}^n E[q_{l} (X^{(l)})]^2 \sum_{k=1}^m \lambda_k p_{kl}
+\sum_{l=1}^n \sum_{j\neq l}^n E[q_{j} (X^{(j)})]E[q_{l} (X^{(l)})] \sum_{k=1}^m \lambda_k p_{kj}p_{kl}.
\end{eqnarray*}

From now on, we assume that $\textbf{q}$ changes with time. Besides, the insurer can control the reserves by paying out dividends with both transaction costs and taxes. That is, there will be a fixed transaction cost $K>0$ and a tax
rate $1-k$ $(0<k<1)$ when the dividends are paid out.
A strategy is described by
$$\alpha=(\textbf{q}_t;\tau_{1},\tau_{2},\ldots, \tau_{n},\ldots; \xi_{1},\xi_{2},\ldots,\xi_{n},\ldots),$$
where $\tau_{n}$ and $\xi_{n}$ denote the times and amounts of dividends.
The controlled surplus process process with strategy $\alpha$ is given by
\begin{eqnarray}\label{e2}
X_t^\alpha=x+\int_0^t [c-\delta(\textbf{q}_s)-a(\textbf{q}_s)] ds+\int_0^t b(\textbf{q}_s) dW_s-\mathop
{\sum}_{n=1}^{\infty}I_{(\tau_{n}<t)}\xi_{n},
\end{eqnarray}
and the corresponding ruin time is defined as
\begin{eqnarray*}
\tau^\alpha=\inf\{t\geq 0: X_t^\alpha< 0\}.
\end{eqnarray*}
\noindent
\begin{definition}\label{Defn:adm-str} \rm {\it
A strategy $\alpha$ is said to be admissible if \begin{enumerate}[(i)]
\item $q_{lt}$ $(l=1,2,\ldots,n)$ are ${\cal F}_{t}$-adapted processes with $0\leq q_{lt}(x) \leq x$ for all $x\geq 0$ and $t\geq 0$.
\item $\tau_{n}$ is a stopping time with respect to  $\{{\cal
F}_{t}\}_{t\geq 0}$ and $ 0\leq \tau_{1}<\tau_{2}<\cdots< \tau_{n}<\cdots a.s.$
\item $\xi_{n}$ is measurable with respect to ${\cal F}_{\tau_{n}-}$ and $0< \xi_{n}\leq X_{\tau_{n}-}^{\alpha}, n=1,2,\ldots$.
\item $P(\lim\nolimits_{n\to \infty}\tau_{n}\leq T)=0$, for all $T\geq 0.$
\end{enumerate} }
\end{definition}

Denoted by $\Pi$ the set of all admissible control strategies.
For a given admissible strategy $\alpha$, we define the performance function (or value function) as
\begin{eqnarray*}\label{e4}
V_\alpha(x)=E\Big[\mathop
{\sum}_{n=1}^{\infty} e^{-\delta \tau_{n}}(k\xi_{n}-K)
I_{\{\tau_{n}<\tau^{\alpha}\}}\mid X_{0-}=x\Big]= E_{x}\Big[\mathop
{\sum}_{n=1}^{\infty} e^{-\delta \tau_{n}}(k\xi_{n}-K)
I_{\{\tau_{n}<\tau^{\alpha}\}}\Big],~
\end{eqnarray*}
which represents the expected total discounted dividends received by the shareholders until the ruin time when the initial surplus is $x$,
where $\delta> 0$ is a priori given discount factor.
Our aim is to find the optimal performance function given by
\begin{eqnarray}\label{e5}
V(x)=\sup_{\alpha\in\Pi}V_\alpha(x),
\end{eqnarray}
and to find the optimal strategy $\alpha^*$ such that $V(x)=V_{\alpha^*}(x)$ for all $x\geq 0$.

In this paper, we assume that the premium is calculated according to the expected value premium principle. For the $l$th ($l=1,2,\ldots,n$) class of insurance business, the positive safety loading for the insurer and reinsurer are $\eta_l$ and $\theta_l$, respectively. Non-cheap reinsurance is considered, that is, $\theta_l>\eta_l$. In order to derive closed-form expressions for the value function $V(x)$ and the corresponding optimal strategy $\alpha^*$, we consider the case of $n=2$ only. Let
\begin{eqnarray*}
 c_l=\sum_{k=1}^m \lambda_k p_{kl}, ~l=1,2, \qquad \mbox{and} \qquad c_3=\sum_{k=1}^m \lambda_k p_{k1}p_{k2}.
\end{eqnarray*}
Then we have
\begin{eqnarray} \label{e6}
&& c=\sum_{l=1}^2 c_l \mu_l (1+\eta_l),~~~
\delta(\textbf{q})=\sum_{l=1}^2 c_l (\mu_l-E[q_{l} (X^{(l)})]) (1+\theta_l), \nonumber \\
&& d(\textbf{q})\triangleq c-\delta(\textbf{q})-a(\textbf{q})=\sum_{l=1}^2 c_l \{\theta_l E[q_{l} (X^{(l)})]-(\theta_l-\eta_l)\mu_l \},\\
&& b^2(\textbf{q})=\sum_{l=1}^2 c_l E[q_{l} (X^{(l)})]^2 +2 c_3 E[q_{1} (X^{(1)})]E[q_{2} (X^{(2)})] .  \nonumber
\end{eqnarray}

\section{The optimal reinsurance form}\label{sec:gain-excess}
\setcounter{equation}{0}
There exists a variety of reinsurance forms in the literature, such as proportional reinsurance, excess-of-loss reinsurance, stop-loss reinsurance, and so on. In this section, we show that the excess-of-loss reinsurance strategy is the optimal reinsurance form for the problem of study.

\begin{lemma}\label{lemma:1}
For any admissible strategy $\alpha=(q_1, q_2;\tau_{1},\ldots, \tau_{n},\ldots; \xi_{1},\ldots,\xi_{n},\ldots)$, there exists an admissible strategy
$\alpha^e=(q_1^e, q_2^e;\tau_{1},\ldots, \tau_{n},\ldots; \xi_{1},\ldots,\xi_{n},\ldots)$ such that $V_\alpha(x)\leq V_{\alpha^e}(x)$, where $(q_1^e, q_2^e)$ is
a two-dimensional excess-of-loss reinsurance strategy.
\end{lemma}
\noindent{\bf Proof.} \ Similar to the proof of Proposition 2.1 of Bai et al. (2013), we know that for any one-dimensional reinsurance strategy $q(\cdot)$ with $0\leq q(x)\leq x$ for $x\geq 0$ and a nonnegative random variable $Z$, there exists an
excess-of-loss reinsurance strategy $q^e(\cdot)=\min\{\cdot, m\}$ with a retention level $0\leq m \leq \infty$ such that
$$E[q^e(Z)]= E[q(Z)],~~E[q^e(Z)]^2\leq E[q(Z)]^2.$$
Then, for any $\textbf{q}=(q_1(\cdot),q_2(\cdot))$, it follows from \eqref{e6} that, there exists $0\leq m_i \leq \infty, i=1,2$ and $\textbf{q}^e=(q_1^e(\cdot),q_2^e(\cdot))
=(\min\{\cdot, m_1\}, \min\{\cdot, m_2\})$ such that
$$d(\textbf{q}^e)=d(\textbf{q}),~~  b^2(\textbf{q}^e)\leq  b^2(\textbf{q}).$$
On the other hand, it is easy to see that both $d(\textbf{q}^e)$ and $b^2(\textbf{q}^e)$ are increasing with respect to $m_1$ and $m_2$, and that $\lim_{(m_1,m_2)\rightarrow (\infty,\infty)} b^2(\textbf{q}^e)\geq  b^2(\textbf{q})$. As a result,
there exists $m'_i \geq m_i, i=1,2$, and $\textbf{q}^{e'}=(\min\{\cdot, m'_1\}, \min\{\cdot, m'_2\})$ such that
$$d(\textbf{q}^{e'})\geq d(\textbf{q}),~~  b^2(\textbf{q}^{e'})=  b^2(\textbf{q}).$$
By \eqref{e2}, we have $X_t^\alpha \leq X_t^{\alpha^{e'}}$. This implies $\tau^\alpha \leq \tau^{\alpha^{e'}}$, which in turn
yields $V_\alpha(x)\leq V_{\alpha^{e'}}(x)$. ~~~~ $\Box$

\vspace{3mm}

Due to Lemma \ref{lemma:1}, we only consider the excess-of-loss reinsurance in the rest of this paper. For notational convenience, we define the following functions:
\begin{eqnarray*}
&& g_l(q)=E(X^{(l)}\wedge q)=\int_0^q \bar{F}_{l}(x) dx,~~ l=1,2,\\
&& G_l(q)=E(X^{(l)}\wedge q)^2=\int_0^q 2x \bar{F}_{l}(x) dx,~~ l=1,2,
\end{eqnarray*}
where $q\in [0,\infty]$ and $\bar{F}_{l}(x)=1-F_{l}(x)=P(X^{(l)}>x)$. Then we have
\begin{eqnarray} \label{e6a}
&& d(\textbf{q})=\sum_{l=1}^2 c_l \{\theta_l g_l(q_l)-(\theta_l-\eta_l)\mu_l \}, \nonumber \\
&& b^2(\textbf{q})=\sum_{l=1}^2 c_l G_l(q_l) +2 c_3 g_1(q_1)g_2(q_2) ,
\end{eqnarray}
where $\textbf{q}=(q_1, q_2)$ with $0\leq q_l \leq \infty$.

\section{QVI and verification theorem}\label{sec:verifi-theo}
\setcounter{equation}{0}
\indent Since the optimal control problem \eqref{e5} is a mixed classical-impulse stochastic control problem,
we deal with it by the method of quasi-variational inequalities (QVI).

For a function $\phi:[0,\infty)\mapsto [0,\infty),$  we
define the maximum operator ${\cal M}$ by
\begin{eqnarray*}\label{M}
{\cal M} \phi(x):=\sup \{\phi(x-\eta)+k\eta -K: 0< \eta \leq x\},
\end{eqnarray*}
and the operator ${\cal L}^{\textbf{q}}$ by
\begin{eqnarray*}
{\cal L}^{\textbf{q}}\phi(x):=\frac{1}{2} b^2(\textbf{q}) \phi''(x)+ d(\textbf{q}) \phi'(x).
 \end{eqnarray*}

Similar to Chen and Yuen (2016), if the value function of (\ref{e5}) is sufficiently smooth, then it satisfies the following QVI:
\begin{eqnarray} \label{e31}
\max \Big\{\max_{0\leq q_1, \, q_2\leq \infty}{\cal L}^{\textbf{q}}V(x)-\delta V(x),~ {\cal M} V(x)-V(x) \Big\}=0,\quad x> 0,
\end{eqnarray}
with boundary condition $V(0)=0$.

\begin{remark} 
Intuitively, the economic insight behind \eqref{e31} reads as: every time when the surplus level $x (>0)$ is such that ${\cal M} V(x) = V (x)$, it is optimal for the insurer to pay lump sum dividends rather than buying reinsurance; while when the current surplus $x (>0)$ is such that ${\cal M} V(x) < V (x)$, it is optimal for the insurer to buy reinsurance to cede out claims rather than paying dividends.
\end{remark}

Furthermore, given a solution $v(x)$ to \eqref{e31}, we can construct the following Markov control strategy.

\begin{definition}\label{Defn:associa-str} \rm {\it
The strategy
$\alpha^{v}=(q_1^{v},q_2^{v};\tau_{1}^{v},\tau_{2}^{v},\cdots, \tau_{n}^{v},\cdots;
\xi_{1}^{v},\xi_{2}^{v},\cdots,\xi_{n}^{v},\cdots )$ is called the
QVI strategy associated with $v$  if the associated process $X^{v}$
given by \eqref{e2} with  $x\geq 0$ satisfies
\begin{eqnarray*}
&& (q_{1t}^{v},q_{2t}^{v})={\rm arg} \mathop {{\max}}_{0\leq q_1, \, q_2\leq \infty} {\cal L}^{\textbf{q}} v(X^{v}_t)~~ \hbox{on} ~~\{v(X_{t}^{v})> {\cal M} v(X_{t}^{v})\}, \\
&& \tau_{1}^{v}={\rm inf}\{t\geq 0: v(X_{t}^{v})= {\cal M} v(X_{t}^{v})\}, \\
&& \xi_{1}^{v}={\rm arg} \mathop {{\rm sup}}_{0<\eta\leq
X^{v}_{\tau_{1}^{v}}}\{v(X^{v}_{\tau_{1}^{v}}-\eta)+k \eta-K\},
\end{eqnarray*}
and for every $n\geq 2$,
\begin{eqnarray*}
\tau_{n}^{v}={\rm inf}\{t> \tau_{n-1}^{v}: v(X_{t}^{v})= {\cal M} v(X_{t}^{v})\},~~\\
\xi_{n}^{v}={\rm arg} \mathop {{\rm sup}}_{0<\eta\leq
X^{v}_{\tau_{n}^{v}}}\{v(X^{v}_{\tau_{n}^{v}}-\eta)+ k\eta-K\}.
\end{eqnarray*}}
\end{definition}

Mimicking the proof of Theorem 3.2 in Chen and Yuen (2016), one can prove the following verification theorem.

\begin{theorem}[Verification Theorem] \label{theo:verifi}
Let $v(x)\in C^1((0,\infty))$ be a solution to \eqref{e31} at all the points with the possible exception of some point where the second derivative may not exist. Suppose there exists $U>0$ such that
$v(x)$ is twice continuously differentiable on $(0,U)$ and $v(x)$ is linear on $[U,\infty)$.
Then $ V(x)\leq v(x),~ x\geq 0.$  Furthermore, if the QVI strategy $\alpha^{v}$ associated
with $v(x)$ is admissible, then $v(x)$ coincides with the
value function $V(x)$ and $\alpha^{v}$ is the optimal strategy, i.e.,
$V(x)=v(x)=V_{\alpha^{v}} (x), \ \ x\geq 0.$
\end{theorem}

\section{Solution to QVI}\label{sec:solution}
\setcounter{equation}{0}

Inspired by Theorem \ref{theo:verifi},
we first assume that there exists a strictly increasing solution $W(x)$ to (\ref{e31}) which is continuously differentiable on $(0,\infty)$ and twice continuously differentiable on $(0,x_1)$,
where $x_1=\inf\{x\geq 0: {\cal M} V(x)= V(x) \}$ (all of these will be proved later).
Then (\ref{e31}) with $V$ replaced by $W$ for $0\leq x< x_1$ can be rewritten as
\begin{eqnarray}\label{e32}
\max_{0\leq q_1, \, q_2\leq \infty} \left\{\frac{1}{2} b^2(\textbf{q}) W''(x)+ d(\textbf{q}) W'(x)-\delta W(x)\right\}=0.
\end{eqnarray}
Let $q_1(x)$ and $q_2(x)$ be the maximizer of the left-hand side of (\ref{e32}).  Assume that $q_1(x)$ and $q_2(x)$ fall in the interval $(0, \infty)$. Differentiating (\ref{e32}) with respect to $q_1$ and $q_2$ respectively, we obtain
\begin{eqnarray}
&& -\frac{W''(x)}{W'(x)}=\frac{~c_1 \theta_1}{c_1 q_1(x)+c_3 g_2[q_2(x)]}, \label{e33} \\
&& -\frac{W''(x)}{W'(x)}=\frac{~c_2 \theta_2}{c_2 q_2(x)+c_3 g_1[q_1(x)]}. ~~~~~\label{e033}
\end{eqnarray}
It follows that
\begin{eqnarray} \label{e331}
\theta_2 q_1(x)-\frac{c_3}{c_2} \theta_1 g_1[q_1(x)]=\theta_1 q_2(x)-\frac{c_3}{c_1} \theta_2 g_2[q_2(x)].
\end{eqnarray}
Let
\begin{eqnarray*}
l_1(q)=\theta_2 q-\frac{c_3}{c_2} \theta_1 g_1(q),~~l_2(q)=\theta_1 q-\frac{c_3}{c_1} \theta_2 g_2(q),~~q\geq 0.
\end{eqnarray*}
Without loss of generality, we assume that $\theta_1\geq \theta_2$. We further assume that
$0=F_l(0)<F_l(x)<1$ for $x>0$ and $l=1,2$. Then it is easy to see that $l_2(q)$ is strictly increasing on $[0,\infty]$,
so the inverse function $l_2^{-1}(q)$ exists. By (\ref{e331}), we have $q_2(x)=l_2^{-1}[l_1(q_1(x))]\geq 0$ if $l_1(q_1(x))\geq 0$.
Let
$$z_l=\sup\{x\geq 0: l_1(x)=0\}.$$
It is easy to see that $0\leq z_l<\infty$ and $l_1(q)\leq 0$ for $q\leq z_l$ since $l''_1(q)\geq 0$. Naturally, we need to find some $q_1(x)\geq z_l$ to guarantee that $q_2(x)\geq 0$ for $x\geq 0$.

Substituting (\ref{e33}) into (\ref{e32}) and replacing $q_2(x)$ with $l_2^{-1}[l_1(q_1(x))]$, we obtain
\begin{eqnarray}\label{e34}
H(q_1(x)) W'(x)-\delta W(x)=0,~~~~
\end{eqnarray}
where
\begin{eqnarray}\label{e3400}
&& H(q)=\sum_{l=1}^2 c_l (\eta_l-\theta_l)\mu_l+c_1 \theta_1 g_1(q)+c_2 \theta_2 g_2[l_2^{-1}l_1(q)] \nonumber\\
&& ~~~~~~~~~~ ~~ -\frac{\theta_1}{2} \frac{c_1 G_1(q)+c_2 G_2[l_2^{-1}l_1(q)]+2 c_3 g_1(q) g_2[l_2^{-1}l_1(q)]}{q+\frac{c_3}{c_1} g_2[l_2^{-1}l_1(q)]}.
\end{eqnarray}
In view of $W(0)=0$ and (\ref{e34}), we see that $H(q_1(0))=0$. So we should discuss the existence of the solution to $H(q)=0$. Now we define an auxiliary following function:
$$k(x)=c_1 \theta_1\Big[g_1(x)-\frac{G_1(x)}{2x} \Big]+k_0,$$
where $k_0=\sum_{l=1}^2 c_l (\eta_l-\theta_l)\mu_l<0$. Since $k'(x)=c_1 \theta_1 \frac{G_1(x)}{2x^2}>0$ for all $x>0$, the inverse function $k^{-1}(x)$ exists. Note that $k(0+)=k_0<0$ and $k(\infty)=c_1 \eta_1\mu_1+c_2 (\eta_2-\theta_2)\mu_2$. Define the zero of $k(x)$ as
\begin{eqnarray*}
z_k=
\left\{
\begin{array}{lll}
k^{-1}(0),~& \theta_2\leq \eta_2+\frac{c_1 \mu_1 \eta_1}{c_2 \mu_2}, \\
~\infty,~& \hbox{otherwise}.
\end{array}
\right.
\end{eqnarray*}

\begin{lemma} \label{lemma1}
There exists a unique solution $q_0$ to $H(q)=0$ on $[z_l,\infty)$ if and only if $z_l\leq z_k$. Furthermore, we have $q_0>0$ if it exists.
\end{lemma}
\noindent{\bf Proof.}~ By some direct calculation, one can show that for $q\geq z_l$,
\begin{eqnarray}\label{e35}
H'(q)&=&\frac{\theta_1}{2} \{c_1 G_1(q)+c_2 G_2[l_2^{-1}l_1(q)]+2 c_3 g_1(q) g_2[l_2^{-1}l_1(q)]\} \nonumber \\
&&  ~~~\times \frac{c_1+c_3 \bar{F}_2[l_2^{-1}l_1(q)] (l_2^{-1}l_1)'(q) }{(c_1 q+ c_3 g_2[l_2^{-1}l_1(q)])^2}.
\end{eqnarray}
On the other hand, for $q> z_l$, we have
\begin{eqnarray*}
 0< l_1(q)=\theta_2 q-\frac{c_3}{c_2} \theta_1 g_1(q)\leq [\theta_2-\frac{c_3}{c_2} \theta_1 \bar{F}_1(q)]q=l'_1(q) q,
\end{eqnarray*}
which implies that $l'_1(q)>0$. As a result, we get $H'(q)>0$ for $q> z_l$, which in turn implies that $H(q)$ is strictly increasing
on $[z_l,\infty]$. Since $l_1(z_l)=0$, we have
\begin{eqnarray*}
H(z_l)&=& k_0+c_1 \Big[\theta_1 g_1(z_l)-\frac{c_2 \theta_2 G_1(z_l)}{2c_3 g_1(z_l)} \Big]\\
&=& k_0+c_1 \theta_1\Big[g_1(z_l)-\frac{G_1(z_l)}{2z_l} \Big]=k(z_l).
\end{eqnarray*}
Besides, we note that
$H(\infty)=\sum_{l=1}^2 c_l \eta_l \mu_l>0$ and $k(x)$ is strictly increasing. It is easy to see that
there exists a unique solution $q_0$ to $H(q)=0$ on $[z_l,\infty)$ if and only if $z_l\leq z_k$.

Furthermore, we have $q_0=z_l$ if $z_l=z_k$ and $q_0>z_l$ if $z_l<z_k$. Note that $z_l\geq 0$ and $z_k>0$. Then we obtain $q_0>0$.
\hfill{$\Box$}

According to Lemma \ref{lemma1}, we will consider the problem in two cases: (1) \ $z_l\leq z_k$;  (2) \ $z_l> z_k$.

\vspace{3mm}

\subsection{The case of $z_l\leq z_k$}

In this case, it follows from Lemma \ref{lemma1} and (\ref{e34}) that $q_1(0)=q_0$. Furthermore,
differentiating (\ref{e34}) with respect to $x$, we have
\begin{eqnarray}\label{e36}
[H'(q_1(x))q'_1(x)-\delta]W'(x)+H(q_1(x))W''(x)=0.
\end{eqnarray}
Using (\ref{e33}) and $q_2(x)=l_2^{-1}[l_1(q_1(x))]$ once again, we obtain
\begin{eqnarray}\label{e37}
W'(x)\Big\{ H'(q_1(x))q'_1(x)- \delta-H(q_1(x)) \frac{c_1 \theta_1}{c_1 q_1(x)+c_3 g_2[l_2^{-1}l_1(q_1(x))]} \Big\}=0.
\end{eqnarray}
Since $W'(x)> 0$, \eqref{e37} gives
\begin{eqnarray}\label{e38}
q'_1(x)=\frac{~\delta+ H(q_1(x)) \frac{c_1 \theta_1}{c_1 q_1(x)+c_3 g_2[l_2^{-1}l_1(q_1(x))]} ~} {H'(q_1(x))}.
\end{eqnarray}
Let
\begin{eqnarray}\label{e39}
G(q)=\int_{q_0}^q \frac{H'(y)} {~\delta+ H(y) \frac{c_1 \theta_1}{c_1 y+c_3 g_2[l_2^{-1}l_1(y)]} ~} dy,~q\geq q_0.
\end{eqnarray}
Since the integrand on the right-hand side of (\ref{e39}) is positive on $[q_0,\infty]$, we see that
$G(q)$ is increasing on $[q_0,\infty]$, and hence the inverse of $G(q)$ exists on $[q_0,\infty]$. As a result, we have
$$q_1(x)=G^{-1}(x), \quad q_2(x)=l_2^{-1}[l_1(G^{-1}(x))]. $$

\begin{lemma} \label{lemma2}
Let $G(q)$ be given by \eqref{e39}. Then we have $G(\infty)<\infty$, which implies that there exists a $x_0=G(\infty)<\infty$ such that
$q_1(x_0)=\infty$.
\end{lemma}
\noindent{\bf Proof.}~ Note that
\begin{eqnarray*}
(l_2^{-1}l_1)'(q)=\frac{1}{\theta_1-\frac{c_3\theta_2}{c_1} \bar{F}_2[l_2^{-1}l_1(q)] }
\times \Big(\theta_2-\frac{c_3\theta_1}{c_2} \bar{F}_1(q)\Big)\rightarrow \frac{\theta_2}{\theta_1},
\end{eqnarray*}
as $q \rightarrow \infty$.  Then it follows from \eqref{e35} that $H'(y)$ tends to 0 at the rate $y^{-2}$ as $y\rightarrow \infty$.
On the other hand, the denominator of the integrand of (\ref{e39}) tends to $\delta$ as $y\rightarrow \infty$.
Then it is easy to see that
$G(\infty)<\infty$, which in turn implies that there exists a $x_0<\infty$ such that $q_1(x_0)=\infty$.
\hfill{$\Box$}

\begin{remark}
Lemma \ref{lemma2} suggests that the insurer will not buy reinsurance when the reserve is no less than $x_0$.
\end{remark}

Assume that $x_0<x_1$ (this will be proved later). Then for $0<x<x_0$, it follows from \eqref{e33} that
\begin{eqnarray} \label{e313}
W(x)=c_4 \int_0^x \exp \Big(- \int_{x_0}^z   \frac{~c_1 \theta_1}{c_1 G^{-1}(y)+c_3 g_2[l_2^{-1}l_1(G^{-1}(y))]} dy \Big) dz, \end{eqnarray}
where $c_4>0$ is a constant.

For $x_0\leq x\leq x_1$, we guess that $q_1(x)=q_2(x)=\infty$. Let
\begin{eqnarray*}
K_1=\frac{1}{2} \sum_{l=1}^2 c_l (\mu_l^2+\sigma_l^2) + c_3 \mu_1 \mu_2, \quad K_2=\sum_{l=1}^2 c_l \eta_l\mu_l.
\end{eqnarray*}
Then (\ref{e32}) becomes
\begin{eqnarray*}
K_1 W''(x)+K_2 W'(x)-\delta W(x)=0,
\end{eqnarray*}
which has the following general solution
\begin{eqnarray}\label{e3131}
W(x)=c_5 e^{r_+(x-x_0)}+c_6 e^{r_-(x-x_0)},
\end{eqnarray}
where $c_5$ and $c_6$ are constants, and
\begin{eqnarray*}
r_+=\frac{-K_2+\sqrt{K_2^2+4\delta K_1}}{2 K_1}, \quad r_-=\frac{-K_2-\sqrt{K_2^2+4\delta K_1}}{2 K_1} .
\end{eqnarray*}

For $x> x_1$, by the definition of $x_1$, we guess that
\begin{eqnarray}\label{e3132}
W(x)=W(\tilde{x})+k(x-\tilde{x})-K,
\end{eqnarray}
where $\tilde{x}<x_1$ is a constant.

By the continuity of $W'$ and $W''$ at $x_0$, it is easy to see that
\begin{eqnarray*}
c_5 r_+ +c_6 r_-=c_4, \quad c_5 r_+^2 +c_6 r_-^2=0,
\end{eqnarray*}
which results in $c_5=c_4 b_1$ and $c_6=c_4 b_2,$ where
\begin{eqnarray}\label{e3133}
b_1=\frac{r_-}{r_+(r_--r_+)} >0, \quad b_2=\frac{r_+}{r_-(r_+-r_-)} <0.
\end{eqnarray}
The unknown constants $c_4$, $\tilde{x}$ and $x_1$ can be determined in the same way as that in Chen and Yuen (2016).  For details, see Chen and Yuen (2016).  The following steps briefly describe how these constants can be determined:
\begin{itemize}
\item[(i)] Define an auxiliary function $U(x)$ as
\begin{eqnarray*}
U(x)=
\left\{
\begin{array}{lll}
\exp\big(- \int_{x_0}^x   \frac{~c_1 \theta_1}{c_1 G^{-1}(y)+c_3 g_2[l_2^{-1}l_1(G^{-1}(y))]} dy\big),~& 0\leq x\leq x_0,  \\
b_1 r_+ e^{r_+(x-x_0)}+b_2 r_- e^{r_-(x-x_0)},~& x> x_0,
\end{array}
\right.
\end{eqnarray*}
which is convex on $(0, \infty)$, and attains its minimum at $x=x_0$ with $U(x_0)=1$.
\item[(ii)] For any fixed $c\in (0,k]$, there exists a unique $\hat{x}_{c}\geq x_0$ such that $cU(\hat{x}_{c})=k$.
Let $\bar{c}=k/U(0)<k$. If $c\in [\bar{c},k]$, then there exists a unique $\tilde{x}_c \in [0,x_0]$ such that $cU(\tilde{x}_c)=k$.
\item[(iii)] Let
\begin{eqnarray*}
I_1(c)=\int_{\tilde{x}_c}^{\hat{x}_{c}} (k-c U(y))dy, \quad c\in [\bar{c}, k],\\
I_2(c)=\int_{0}^{\hat{x}_{c}} (k-c U(y))dy, \quad c\in [0, k].
\end{eqnarray*}
If $I_1(\bar{c})>K$, then there exists a unique $c^* \in (\bar{c}, k)$ such that $I_1(c^*)=K$.
If $I_1(\bar{c})\leq K$, then there exists a unique $c^* \in (0, k)$ such that $I_2(c^*)=K$.
\item[(iv)] Let $c_4=c^*$, $x_1=\hat{x}_{c^*}>x_0$, and $\tilde{x}=\tilde{x}_{c^*}$, where
$\tilde{x}_{c^*}=0$ if $I_1(\bar{c})\leq K$.
\end{itemize}

These together with \eqref{e313}-\eqref{e3132} yield
\begin{equation}\label{e314}
W(x)=
\left\{
\begin{array}{lll}
c^* \int_0^x \exp \big(- \int_{x_0}^z   \frac{~c_1 \theta_1}{c_1 G^{-1}(x)+c_3 g_2[l_2^{-1}l_1(G^{-1}(x))]} dy \big) dz,~& 0\leq x< x_0,\\
c^* [b_1 e^{r_+(x-x_0)}+b_2 e^{r_-(x-x_0)}],~ & x_0\leq x< \hat{x}_{c^*},\\
W(\tilde{x}_{c^*})+k(x-\tilde{x}_{c^*})-K, ~ & x\geq \hat{x}_{c^*},
\end{array}
\right.
\end{equation}
where $b_1$ and $b_2$ are given in \eqref{e3133}.

\begin{theorem} \label{theo:impul-1}
If $z_l\leq z_k$, then the function $W(x)$ of \eqref{e314} is continuously differentiable on $(0, \infty)$
and twice continuously differentiable on $(0, \hat{x}_{c^*})\cup (\hat{x}_{c^*}, \infty)$.
Furthermore, $W(x)$ is a solution to the QVI of \eqref{e31}.
\end{theorem}
\noindent{\bf Proof.}~ One can prove the theorem by replacing $G(1)$ and $\max_{0\leq b\leq 1,\, 0\leq u\leq 1}{\cal L}^{b,u}W(x)$ by $x_0$ and $\max_{0\leq q_1, \, q_2\leq \infty}{\cal L}^{\textbf{q}}W(x)$, respectively, and then mimicking the steps in the proof of Theorem 4.1 of Chen and Yuen (2016).  \hfill{$\Box$}

\vspace{3mm}

\subsection{The case of $z_l> z_k$}

In this case, it follows from Lemma \ref{lemma1} that the equation $H(q)=0$ on $[z_l,\infty)$ has no solution. Then we guess that
$q_2(x)=0$. Then (\ref{e32}) becomes
\begin{eqnarray}\label{e521}
\max_{0\leq q_1 \leq \infty} \left\{\frac{1}{2} c_1 G_1(q_1) W''(x)+ [c_1 \theta_1 g_1(q_1)+k_0] W'(x)-\delta W(x)\right\}=0.
\end{eqnarray}
Differentiating (\ref{e521}) with respect to $q_1$, we obtain
\begin{eqnarray*}
c_1 \bar{F}_1(q_1)[q_1 W''(x)+\theta_1 W'(x)]=0,
\end{eqnarray*}
which yields
\begin{eqnarray} \label{e522}
\frac{W''(x)}{W'(x)}=-\frac{~ \theta_1}{q_1(x)}.
\end{eqnarray}
Substituting (\ref{e522}) into (\ref{e521}), we obtain
\begin{eqnarray}\label{e523}
k(q_1(x)) W'(x)-\delta W(x)=0.
\end{eqnarray}
Differentiating (\ref{e523}) with respect to $x$ and using (\ref{e522}) once again, we have
\begin{eqnarray}\label{e57}
W'(x)\Big\{ k'(q_1(x))q'_1(x)- \delta- \frac{\theta_1 k(q_1(x))}{ q_1(x)} \Big\}=0.
\end{eqnarray}
Since $W'(x)> 0$, \eqref{e57} gives
\begin{eqnarray*}
q'_1(x)=\frac{~\delta+ \frac{\theta_1 k(q_1(x))}{ q_1(x)} ~} {k'(q_1(x))}.
\end{eqnarray*}
In view of $W(0)=0$ and (\ref{e523}), we see that $k(q_1(0))=0$, which implies that $q_1(0)=z_k>0$.
Let
\begin{eqnarray}\label{e59}
R_1(q)=\int_{z_k}^q \frac{k'(y)} {~\delta+  \frac{\theta_1 k(y)}{y } ~} dy,~q\geq z_k.
\end{eqnarray}
Since the integrand on the right-hand side of (\ref{e59}) is positive on $[z_k,\infty]$, $R_1(q)$ is increasing on $[z_k,\infty]$, which implies that the inverse of $R_1(q)$ exists on $[z_k,\infty]$. Let $\tilde{x}_0=R_1(z_l)$. Then for $0<x\leq \tilde{x}_0$, we have $q_1(x)=R_1^{-1}(x),~ q_2(x)=0$, and it follows from (\ref{e522}) that
\begin{eqnarray} \label{e510}
W(x)=C_1 \int_0^x \exp \Big(- \int_{\tilde{x}_0}^z   \frac{~\theta_1}{R_1^{-1}(y)} dy \Big) dz, \end{eqnarray}
where the constant $C_1>0$ will be determined later.

For $x>\tilde{x}_0$, similar to the case of $z_l\leq z_k$, it can be shown that $q_1(x)$ satisfies \eqref{e38}.
Note that $q_1(\tilde{x}_0)=z_l$. Define
\begin{eqnarray*}
R_2(q)=\int_{z_l}^q \frac{H'(y)} {~\delta+ H(y) \frac{c_1 \theta_1}{c_1 y+c_3 g_2[l_2^{-1}l_1(y)]} ~} dy,~q\geq z_l.
\end{eqnarray*}
Let
$$q_1(x)=R_2^{-1}(x-\tilde{x}_0),  \quad  x>\tilde{x}_0.$$
Similar to Lemma \ref{lemma2}, there exists a $x_0\in (\tilde{x}_0, \infty)$ such that $q_1(x_0)=\infty$.
Then for $\tilde{x}_0<x<x_0$, we have
\begin{eqnarray} \label{e511}
& q_1(x)=R_2^{-1}(x-\tilde{x}_0), \quad q_2(x)=l_2^{-1}[l_1(R_2^{-1}(x-\tilde{x}_0))], \nonumber\\
& W(x)=C_2 \int_{\tilde{x}_0}^x \exp \big(- \int_{x_0}^z   \frac{~c_1 \theta_1}{c_1 R_2^{-1}(y-\tilde{x}_0)+c_3 g_2[l_2^{-1}l_1(R_2^{-1}(y-\tilde{x}_0))]} dy \big)dz +C_3, \end{eqnarray}
where the constants $C_2$ and $C_3 >0$ will be determined later.

For $x\geq x_0$, we guess that $q_1(x)=q_2(x)=\infty$, and $W(x)$ is the same as \eqref{e3131} and \eqref{e3132} for $x_0\leq x\leq x_1$ and $x>x_1$, respectively.

As a result, we have
\begin{equation}\label{e512}
 W(x)=
\left\{
\begin{array}{lll}
C_1 \int_0^x \exp \big(- \int_{\tilde{x}_0}^z   \frac{~\theta_1}{R_1^{-1}(y)} dy \big) dz, \hspace{-0.2cm}& 0\leq x< \tilde{x}_0,\\
C_2 \int_{\tilde{x}_0}^x \exp \big(- \int_{x_0}^z   \frac{~c_1 \theta_1}{c_1 R_2^{-1}(y-\tilde{x}_0)+c_3 g_2[l_2^{-1}l_1(R_2^{-1}(y-\tilde{x}_0))]} dy \big)dz +C_3, \hspace{-0.2cm}& \tilde{x}_0\leq x< x_0,\\
c_5 e^{r_+(x-x_0)}+c_6 e^{r_-(x-x_0)},  \hspace{-0.2cm}& x_0\leq x< x_1,\\
W(\tilde{x})+k(x-\tilde{x})-K, \hspace{-0.2cm}& x\geq x_1.
\end{array}
\right.
\end{equation}

We now need to determine the unknown constants mentioned above.
By the continuity of $W'$ at $\tilde{x}_0$, we have
\begin{eqnarray*}
C_1=C_2  \exp \Big(\int_{\tilde{x}_0}^{x_0} \frac{~c_1 \theta_1}{c_1 R_2^{-1}(y-\tilde{x}_0)+c_3 g_2[l_2^{-1}l_1(R_2^{-1}(y-\tilde{x}_0))]} dy \Big).
\end{eqnarray*}
Besides, \eqref{e511} and (\ref{e523}) imply that $C_3=W(\tilde{x}_0)=\frac{k(z_l)}{\delta}C_1$. Moreover,
the continuity of $W'$ and $W''$ at $x_0$ implies that
$c_5=C_2 b_1$ and $c_6=C_2 b_2,$ where $b_1$ and $b_2$ are given in \eqref{e3133}.
Then it is enough to determine the constants $C_2$, $\tilde{x}$ and $x_1$, which can be obtained by using steps similar to those presented in Section 5.1. Analogous to Theorem \ref{theo:impul-1}, we have the following result.

\begin{theorem} \label{theo:impul-2}
If $z_l> z_k$, then the function $W(x)$ of \eqref{e512} is continuously differentiable on $(0, \infty)$
and twice continuously differentiable on $(0, \hat{x}_{c^*})\cup (\hat{x}_{c^*}, \infty)$.
Furthermore, $W(x)$ is a solution to the QVI of \eqref{e31}.
\end{theorem}
\noindent{\bf Proof.}~ For $x\geq \tilde{x}_0$, the proof is similar to that of Theorem \ref{theo:impul-1}. So, we only prove that $W(x)$ is a solution to \eqref{e31} for $0\leq x< \tilde{x}_0$. Since one can show that $W(x)$ of \eqref{e512} satisfies ${\cal L}^{\textbf{q}^*}W(x)-\delta W(x)=0$ with $\textbf{q}^*=(q_1^*(x), q_2^*(x))=(R_1^{-1}(x), 0)$. As a consequence, we need to show that ${\cal L}^{\textbf{q}}W(x)-\delta W(x)\leq 0$ for any $q_1, q_2\in [0,\infty]$, which is equivalent
to verify that ${\cal L}^{\textbf{q}}W(x)-{\cal L}^{\textbf{q}^*}W(x)\leq 0$ for any $q_1, q_2\in [0,\infty]$.
By (\ref{e522}), the latter is then equivalent to
\begin{eqnarray*}
[d(\textbf{q})-d(\textbf{q}^*)]-\frac{\theta_1}{2q_1^*} [b^2(\textbf{q})-b^2(\textbf{q}^*)]\leq 0.
\end{eqnarray*}
Let $$\varphi(\textbf{q})=d(\textbf{q})-\frac{\theta_1}{2q_1^*} b^2(\textbf{q}),~~\textbf{q} \in [0,\infty]\times [0,\infty].$$
Then it is enough to show that $\varphi(\textbf{q})$ attains its maximum at $\textbf{q}=\textbf{q}^*=(q_1^*, 0)$.
Note that
\begin{eqnarray} \label{e513}
&& \frac{\partial \varphi(\textbf{q})}{\partial q_1}=\frac{c_1 \theta_1 \bar{F}_1(q_1)}{q_1^*} \Big[q_1^*-q_1-\frac{c_3}{c_1} g_2(q_2)\Big],  \nonumber\\
&& \frac{\partial \varphi(\textbf{q})}{\partial q_2}=\frac{\bar{F}_2(q_2)}{q_1^*}
\big[c_2 \theta_2 q_1^*-c_2 \theta_1 q_2-c_3 \theta_1 g_1(q_1)\big].
\end{eqnarray}
Since $\frac{\partial \varphi(\textbf{q})}{\partial q_1}<0$ for $q_1>q_1^*$, we only consider the case of $q_1\leq q_1^*$.
If $\frac{\partial \varphi(\textbf{q})}{\partial q_1}=0$, then $q_2=g_2^{-1}[\frac{c_1}{c_3}(q_1^*-q_1)]$, and
\eqref{e513} yields
$\frac{\partial \varphi(\textbf{q})}{\partial q_2}=\frac{\bar{F}_2(q_2)}{q_1^*} L(q_1)$, where
\begin{eqnarray*}
L(q_1)=c_2 \theta_2 q_1^*-c_3 \theta_1 g_1(q_1)-c_2 \theta_1 g_2^{-1}\Big[\frac{c_1}{c_3}(q_1^*-q_1)\Big],~~0\leq q_1\leq q_1^*.
\end{eqnarray*}
It is easy to see that
\begin{eqnarray*}
L'(q_1)&=& -c_3 \theta_1 \bar{F}_1(q_1)+c_2 \theta_1 \frac{1}{\bar{F}_2\big[\frac{c_1}{c_3}(q_1^*-q_1)\big]}\frac{c_1}{c_3}\\
&=& \frac{\theta_1}{c_3 \bar{F}_2\big[\frac{c_1}{c_3}(q_1^*-q_1)\big]} \Big\{c_1 c_2-c_3^2 \bar{F}_1(q_1)\bar{F}_2\Big[\frac{c_1}{c_3}(q_1^*-q_1)\Big]   \Big\}>0.
\end{eqnarray*}
On the other hand, for $0\leq x< \tilde{x}_0$, we have $z_k\leq q_1^*<z_l$, and
$$L(q_1^*)=c_2 l_1(q_1^*)\leq c_2 l_1(z_l)=0.$$
Therefore, we obtain
\begin{eqnarray*}
\frac{\partial \varphi(\textbf{q})}{\partial q_2}\Big |_{(q_1, q_2)\in \big\{(q_1, q_2):\frac{\partial \varphi(\textbf{q})}{\partial q_1}=0 \big\}}\leq 0.
\end{eqnarray*}
As a result, we see that $\varphi(\textbf{q})$ attains its maximum at $(q_1^*, 0)$.   \hfill{$\Box$}

\section{The value function and the optimal policy}

Let
\begin{equation*}
(q_{1t}^*, q_{2t}^*)=
\left\{
\begin{array}{lll}
\big(G^{-1}(X_t^*), l_2^{-1}[l_1(G^{-1}(X_t^*))]\big),~& 0\leq X_t^*\leq x_0,  \\
(\infty, \infty),~& X_t^*>x_0,
\end{array}
\right.
\end{equation*}
when $z_l\leq z_k$; and
\begin{eqnarray*}
(q_{1t}^*, q_{2t}^*)=
\left\{
\begin{array}{lll}
(R_1^{-1}(X_t^*), 0),~& 0\leq X_t^*< \tilde{x}_0,\\
\big(R_2^{-1}(X_t^*-\tilde{x}_0), l_2^{-1}[l_1(R_2^{-1}(X_t^*-\tilde{x}_0))]\big), ~& \tilde{x}_0\leq X_t^*< x_0,\\
(\infty, \infty),~ & X_t^*\geq x_0,
\end{array}
\right.
\end{eqnarray*}
when $z_l> z_k$.  Recall $X_t^*$ of \eqref{e2} with $\alpha=\alpha^*=(q_{1t}^*, q_{2t}^*;\tau_{1}^*,\tau_{2}^*,\cdots ;
\xi_{1}^*,\xi_{2}^*,\cdots)$. Define $\{\tau_n^*,~ \xi_n^*,~ n\geq 1\}$ as follows:
\begin{itemize}
\item[(i)] If $I_1(\bar{c})>K$, then we define
$$\tau_1^*=\inf\{t> 0: X_t^*= \hat{x}_{c^*} \}, \quad \xi_1^*=\hat{x}_{c^*}-\tilde{x}_{c^*},$$
when the initial surplus $0<x<\hat{x}_{c^*}$,
$$\tau_1^*=0, \quad \xi_1^*=x-\tilde{x}_{c^*},$$
when the initial surplus $x\geq \hat{x}_{c^*}$, and
$$\tau_n^*=\inf\{t> \tau_{n-1}^*: X_t^*= \hat{x}_{c^*} \},~ \xi_n^*=\hat{x}_{c^*}-\tilde{x}_{c^*},$$
for every $n\geq 2$;
\item[(ii)] If $I_1(\bar{c})\leq K$, then we define
$$\tau_1^*=\inf\{t> 0: X_t^*= \hat{x}_{c^*} \}, \quad \xi_1^*=\hat{x}_{c^*},$$
when the initial surplus $0<x<\hat{x}_{c^*}$,
$$\tau_1^*=0, \quad \xi_1^*=x,$$
when the initial surplus $x\geq \hat{x}_{c^*}$, and
$$\tau_n^*=\infty, \quad \xi_n^*=0,$$
for every $n\geq 2$.
\end{itemize}

\begin{theorem} \label{theo:impul-3}
The value function $V(x)$ is given by \eqref{e314} when $z_l\leq z_k$, and by \eqref{e512} when $z_l> z_k$; and the strategy $\alpha^*$ is the corresponding optimal policy.
\end{theorem}
\noindent{\bf Proof.} \ The proof is similar to that of Theorem 4.2 of Chen and Yuen (2016).

\section{Numerical examples}
In this section, we give some numerical examples to assess the impact of some parameters on the optimal reinsurance policy.
We assume that the claim sizes $X^{(1)}$ and $X^{(2)}$ are exponentially distributed with parameters $\beta_1$ and $\beta_2,$ respectively. Then, for $l=1,2,$ we have $\mu_l=\frac{1}{\beta_l}, \sigma_l^2=\frac{1}{\beta_l^2},
g_l(q)=\frac{1}{\beta_l}(1-e^{-\beta_l q})$, and $G_l(q)=\frac{2}{\beta_l^2}[1-(1+\beta_l q)e^{-\beta_l q}]$.
We take $m=3, n=2,$ $p_{11}=p_{22}=1$, $p_{12}=p_{21}=0$ and $p_{31}=p_{32}=1$ so that the resulting model reduces to the common shock model. Besides, we set $\beta_1=1, \beta_2=2, \eta_1=1, \eta_2=0.8, \delta=0.5, \lambda_1=3, \lambda_2=4, \theta_2=1.$
For $\theta_1=1.2$, the effect of $\lambda_3$ on the optimal reinsurance policy is studied in Example \ref{exmple-1}. Example \ref{exmple-2} shows the effect of $\theta_1$ on the optimal reinsurance policy for $\lambda_3=2$.

\begin{example} \label{exmple-1}\rm   In this example, we set $\theta_1=1.2$, and take $\lambda_3=1, 1.5, 2,$ respectively. The effect of $\lambda_3$ on the optimal reinsurance strategies $q_1(x)$ and $q_2(x)$ are shown in Figures \ref{fig1}-\ref{fig3}.

Table \ref{table-1} shows that the critical point $x_0$ increases as $\lambda_3$ increases.
We see from Figures \ref{fig1} and \ref{fig2} that both $q_1(x)$ and $q_2(x)$ are strictly increasing functions, and they decrease as $\lambda_3$ increases. This means that the optimal retention level is higher for larger reserve, and is lower when the insurers face higher risk. The result coincides with our intuition. We also observe from Figure \ref{fig3} that the difference of two reinsurance strategies is quite small, and both of the reinsurance strategies change slowly for small reserve, while they are quite sensitive to the change of surplus when the surplus near the critical point $x_0$.

\begin{table}[!htbp]
\begin{center}
\begin{tabular}{|c|c|c|c|} \hline
$\lambda_3$  & 1  & 1.5   &  2  \\
\hline
$x_0$ & 2.2170 & 2.4666 & 2.7262 \\
\hline
\end{tabular}
\caption{Effect of $\lambda_3$ on the critical point $x_0$.}\label{table-1}
 \end{center}
\end{table}

\begin{figure}[!htbp]
\begin{center}
\includegraphics[width=10cm, height=7cm]{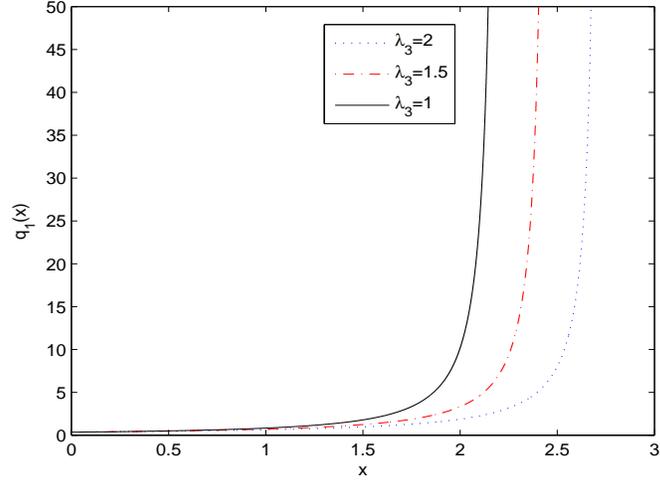}
\caption{Impact of $\lambda_3$ on the optimal reinsurance policy $q_1(x)$.}\label{fig1}
\end{center}
\end{figure}

\begin{figure}[!htbp]
\begin{center}
\includegraphics[width=10cm, height=7cm]{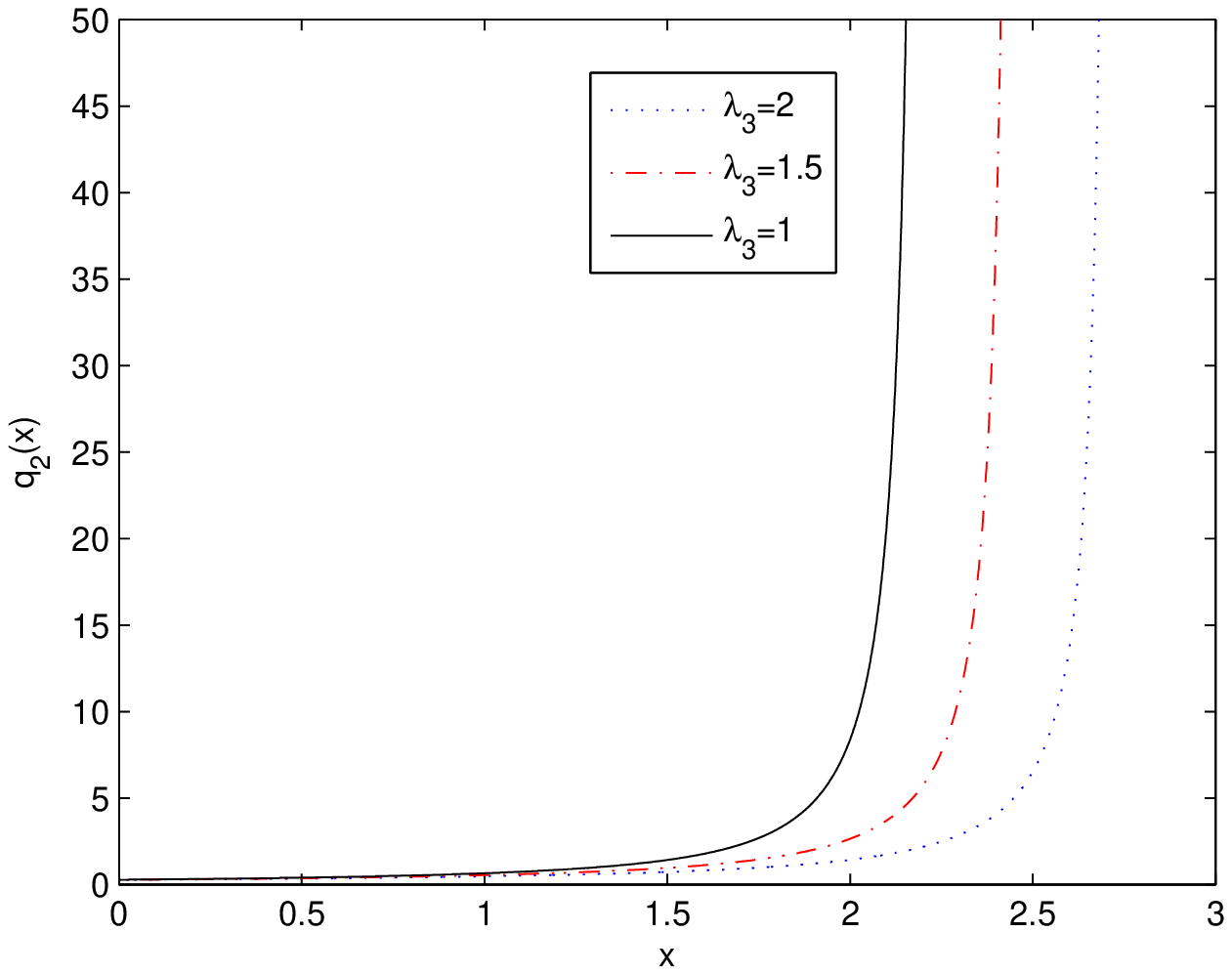}
\caption{Impact of $\lambda_3$ on the optimal reinsurance policy $q_2(x)$.}\label{fig2}
\end{center}
\end{figure}

\begin{figure}[!htbp]
\begin{center}
\includegraphics[width=10cm, height=7cm]{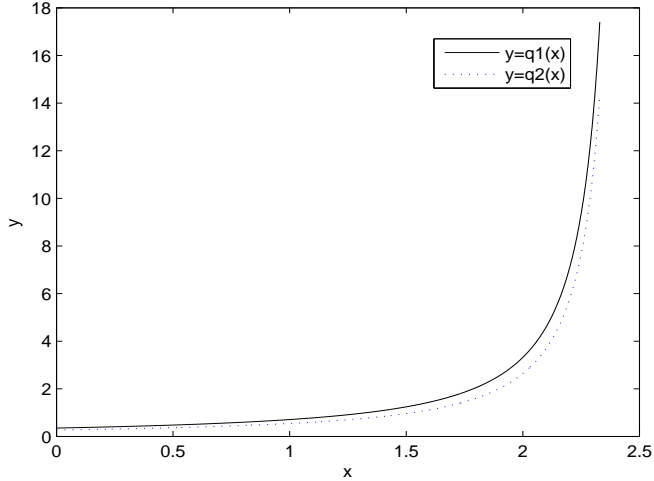}
\caption{The difference of $q_1(x)$ and $q_2(x)$ for $\lambda_3=1.5$.}\label{fig3}
\end{center}
\end{figure}

\end{example}

\begin{example} \label{exmple-2}\rm   In this example, we set $\lambda_3=2$, and take $\theta_1=1.2, 1.5, 2.1,$ respectively.

We see from Table \ref{table-2} that the critical point $x_0$ also increases as $\theta_1$ increases, which means that the insurer should hold a larger reserve when the reinsurance premium becomes more expensive.  Figures \ref{fig4} and \ref{fig5} indicate that both $q_1(x)$ and $q_2(x)$ are not strictly decreasing with respect to $\theta_1$. We can also see that the change of $\theta_1$ has larger effect on $q_1(x)$ than that on $q_2(x)$.  When the reinsurance premium is more expensive, the insurer with small reserve tends to buy less reinsurance, and vice versa.

\begin{table}[!htbp]
\begin{center}
\begin{tabular}{|c|c|c|c|} \hline
$\theta_1$  & 1.2  & 1.5   &  2.1  \\
\hline
$x_0$ & 2.7262 & 4.8197 &  7.8058 \\
\hline
\end{tabular}
\caption{Effect of $\theta_1$ on the critical point $x_0$.}\label{table-2}
 \end{center}
\end{table}

\begin{figure}[!htbp]
\begin{center}
\includegraphics[width=10cm, height=7cm]{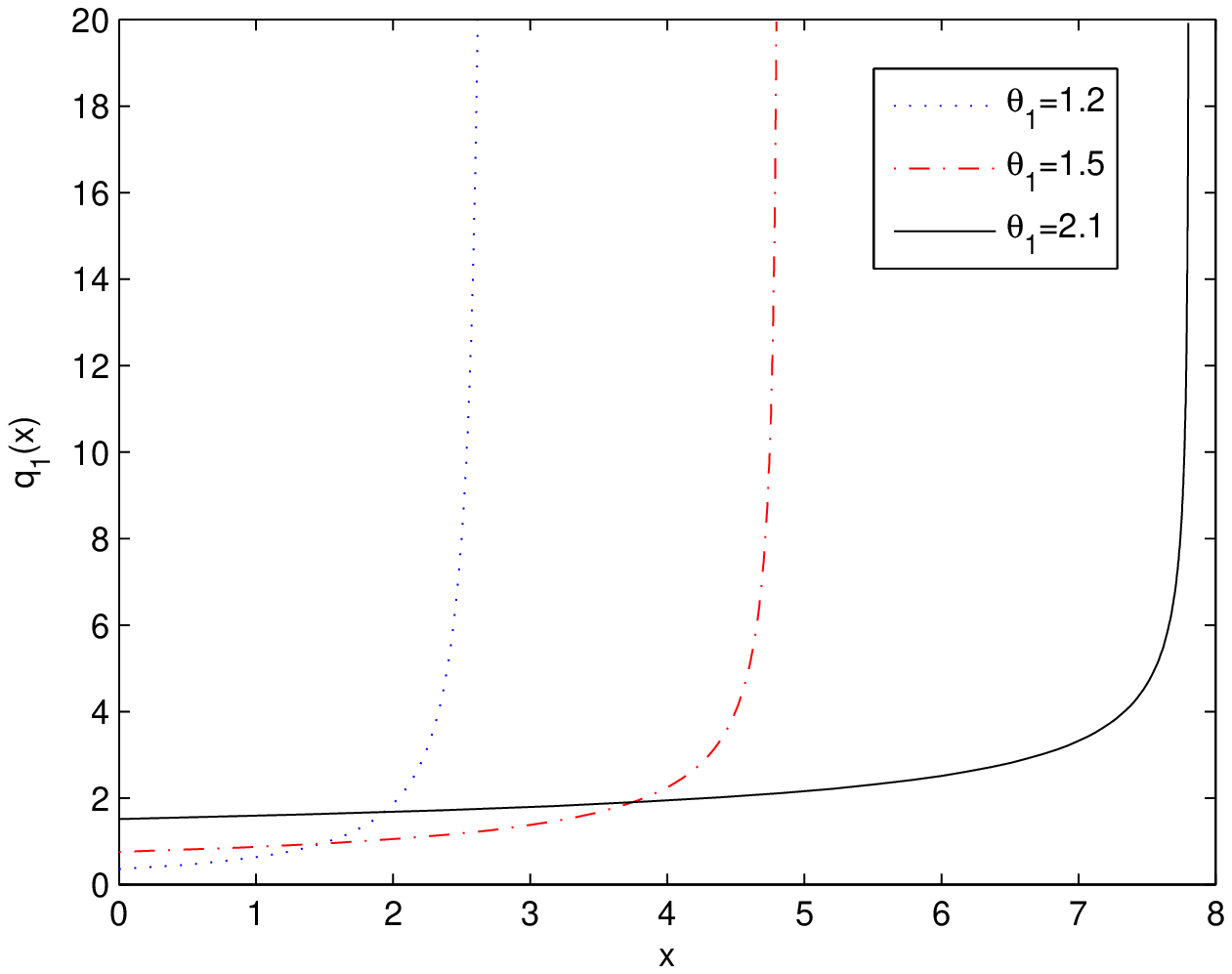}
\caption{Impact of $\theta_1$ on the optimal reinsurance policy $q_1(x)$.}\label{fig4}
\end{center}
\end{figure}

\begin{figure}[!htbp]
\begin{center}
\includegraphics[width=10cm, height=7cm]{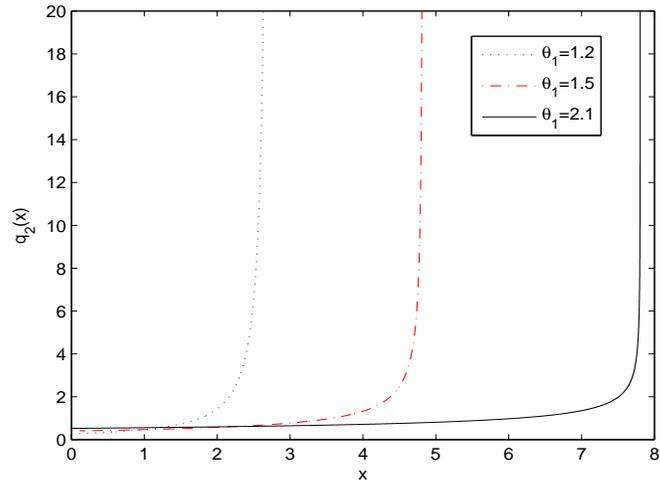}
\caption{Impact of $\theta_1$ on the optimal reinsurance policy $q_2(x)$.}\label{fig5}
\end{center}
\end{figure}

\end{example}

\vspace{0.5 cm}

\noindent {\large \bf Acknowledgements}

\vspace{0.3 cm}

The research of Mi Chen was supported by National Natural Science Foundation of China (Nos. 11701087 and 11701088), Natural Science Foundation of Fujian Province (Nos. 2018J05003 and 2019J01673), Program for Innovative Research Team in Science and
Technology in Fujian Province University, and the grant ``Probability and Statistics: Theory and Application (No. IRTL1704)" from Fujian Normal University. The research of Kam Chuen Yuen was supported by a grant from the Research Grants Council of the Hong Kong Special Administrative Region, China (Project No. HKU17329216). The research of Wenyuan Wang was supported by National Natural Science Foundation of China (No. 11661074).



\vspace{0.5cm}

\noindent{\large \bf References}

\vspace{0.3 cm}

\noindent   Albrecher, H., Thonhauser, S. (2008). Optimal dividend strategies for a risk process under force of interest.
Insurance: Mathematics and Economics 43(1), 134-149. 

\vspace{0.3 cm}

\noindent Albrecher, H., Thonhauser, S. (2009). Optimality results for dividend problems in insurance.
Revista de la Real Academia de Ciencias Exactas, Fisicas y Naturales-Serie A: Matematicas 103(2), 295-320.

\vspace{0.3 cm}

\noindent Asmussen, S., H{\o}jgaard, B. and Taksar, M. (2000).  Optimal risk control and dividend distribution policies. Example of excess-of loss reinsurance. Finance and Stochastics 4(3), 299-324.

\vspace{0.3 cm}

\noindent  Asmussen, S. and Taksar, M. (1997). Controlled diffusion models for optimal dividend pay-out. Insurance: Mathematics and Economics 20(1), 1-15.

\vspace{0.3 cm}

\noindent  Avram, F., Palmowski, Z. and Pistorius, M. (2007). On the optimal dividend problem for a spectrally negative
L\'{e}vy process. The Annals of Applied Probability 17(1), 156-180.

\vspace{0.3 cm}

\noindent Avram, F., Palmowski, Z. and Pistorius, M. (2015). On Gerber-Shiu functions and optimal dividend distribution
for a L\'{e}vy risk process in the presence of a penalty function. The Annals of Applied Probability 25(4), 1868-1935.

\vspace{0.3 cm}

\noindent  Azcue, P. and Muler, N. (2005). Optimal reinsurance and dividend distribution policies in the Cram\'{e}r-Lundberg model. Mathematical Finance 15(2), 261-308.

\vspace{0.3 cm}

\noindent Azcue, P. and Muler, N. (2012). Optimal dividend policies for compound Poisson processes: The case of
bounded dividend rates. Insurance: Mathematics and Economics 51(1), 26-42.

\vspace{0.3 cm}

\noindent Bai, L., Cai, J. and Zhou, M. (2013).  Optimal reinsurance policies for an insurer with a bivariate reserve risk process in a dynamic setting. Insurance: Mathematics and Economics 53(3), 664-670.

\vspace{0.3 cm}

\noindent Bai, L., Guo, J. and Zhang, H. (2010).  Optimal excess-of-loss reinsurance and dividend payments with both transaction costs and taxes.  Quantitative Finance 10(10), 1163-1172.

\vspace{0.3 cm}

\noindent Chen M., Peng, X. and Guo, J. (2013) Optimal dividend problem with a nonlinear regular-singular stochastic control. Insurance: Mathematics and Economics 52(3), 448-456.

\vspace{0.3 cm}

\noindent Chen, M. and Yuen, K.C. (2016). Optimal dividend and reinsurance in the presence of two reinsurers. Journal of Applied Probability 53(2), 554-571.

\vspace{0.3 cm}

\noindent  Czarna, I. and Palmowski, Z. (2010). Dividend problem with parisian delay for a spectrally negative L\'{e}vy risk process. Journal of Optimization Theory and Applications 161(1), 239-256.

\vspace{0.3 cm}

\noindent  Gerber, H.U. and Shiu, E.S.W. (2004). Optimal dividends: analysis with Brownian motion. North American Actuarial Journal 8(1), 1-20.

\vspace{0.3 cm}

\noindent Gerber, H.U. and Shiu, E.S.W. (2006). On optimal dividend strategies in the compound Poisson model. North American Actuarial Journal 10(2), 76-93.

\vspace{0.3 cm}

\noindent Han, X., Liang, Z. and Yuen, K.C. (2018).  Optimal proportional reinsurance to minimize the probability of drawdown under thinning-dependence structure. Scandinavian Actuarial Journal 10, 863-889.

\vspace{0.3 cm}

\noindent He, L. and Liang, Z. (2008).  Optimal financing and dividend control of the insurance company with proportional reinsurance policy. Insurance: Mathematics and Economics 42(3), 976-983.

\vspace{0.3 cm}

\noindent Hernandez, C. and Junca, M. (2015). Optimal dividend payments under a time of ruin constraint: Exponential claims.
Insurance: Mathematics and Economics 65(1), 136-142.

\vspace{0.3 cm}

\noindent  H{\o}jgaard, B. (2002). Optimal dynamic premium control in non-life insurance. Maximizing dividend payouts. Scandinavian Actuarial Journal 4, 225-245.

\vspace{0.3 cm}

\noindent H{\o}jgaard, B. and Taksar, M. (1999).  Controlling risk exposure and dividends payout schemes: insurance company example. Mathematical Finance 9(2), 153-182.

\vspace{0.3 cm}

\noindent Hunting, M. and Paulsen, J. (2013). Optimal dividend policies with transaction costs for a class of jump-diffusion processes. Finance and Stochastics 17(1), 73-106.

\vspace{0.3 cm}

\noindent Jeanblanc-Picqu{\'e}, M. and Shiryaev, A. (1995).  Optimization of the flow of dividends. Russian Mathematical Surveys 50(2), 57-77.

\vspace{0.3 cm}

\noindent  Kyprianou, A. and Palmowski, Z. (2007). Distributional study of de Finetti's dividend problem for a general L\'{e}vy insurance risk process. Journal of Applied Probability 44(2), 428-443.

\vspace{0.3 cm}

\noindent Li, Y., Bai, L. and Guo, J. (2016). Optimal dividend and reinsurance problem for an insurance company with dependent risks (in Chinese).  Scientia Sinica Mathematica 46(8), 1161-1178.

\vspace{0.3 cm}

\noindent Li, Y., Li, Z., Wang, S. and Xu, Z. (2020). Dividend optimization for jump-diffusion model with solvency constraints. Operations Research Letters 48(2), 170-175.

\vspace{0.3 cm}

\noindent Liang, X. and Palmowski, Z. (2018). A note on optimal expected utility of dividend payments with proportional reinsurance. Scandinavian Actuarial Journal 4, 275-293.

\vspace{0.3 cm}

\noindent Liang, Z. and Yuen, K.C. (2016).  Optimal dynamic reinsurance with dependent risks: variance premium principle. Scandinavian Actuarial Journal 1, 18-36.

\vspace{0.3 cm}

\noindent  Loeffen, R. (2008). On optimality of the barrier strategy in de Finetti's dividend problem for spectrally negative L\'{e}vy processes. Journal of Applied Probability 18(5), 1669-1680.

\vspace{0.3 cm}

\noindent  Loeffen, R. (2009). An optimal dividends problem with transaction costs for spectrally negative L\'{e}vy processes. Insurance: Mathematics and Economics 45(1), 41-48.

\vspace{0.3 cm}

\noindent  Loeffen, R. and Renaud, J.F. (2010). De Finetti's optimal dividends problem with an affine penalty function at
ruin. Insurance: Mathematics and Economics 46(1), 98-108.

\vspace{0.3 cm}

\noindent L{\o}kka, A. and Zervos, M. (2008). Optimal dividend and issuance of equity policies in the presence of proportional costs. Insurance: Mathematics and Economics 42(3), 954-961.

\vspace{0.3 cm}

\noindent  Paulsen, J. (2003). Optimal dividend payouts for diffusions with solvency constraints. Finance and Stochastics 7(4), 457-473.

\vspace{0.3 cm}

\noindent  Paulsen, J. (2007). Optimal dividend payments until ruin of diffusion processes when payments are subject to both fixed and proportional costs. Advances in Applied Probability 39(3), 669-689.

\vspace{0.3 cm}

\noindent  Paulsen, J. (2008). Optimal dividend payments and reinvestments of diffusion processes with both fixed and proportional costs. SIAM Journal on Control \& Optimization 47(5), 2201-2226.

\vspace{0.3 cm}

\noindent Peng, X., Bai, L. and Guo, J. (2016). Optimal control with restrictions for a diffusion risk model under constant interest force. Applied Mathematics \& Optimization 73(1), 115-136.

\vspace{0.3 cm}

\noindent P\'{e}rez, J., Yamazaki, K. and Yu, X. (2018). On the bail-out optimal dividend problem.
Journal of Optimization Theory and Applications 179(2), 553-568.

\vspace{0.3 cm}

\noindent  Schmidli, H. (2006). Optimisation in non-life insurance. Stochastic Models 22(4), 689-722. 

\vspace{0.3 cm}

\noindent Vierk\"{o}tter, M. and Schmidli, H. (2017). On optimal dividends with exponential and linear penalty payments.
Insurance: Mathematics and Economics 72(1), 265-270.

\vspace{0.3 cm}

\noindent Wang, G. and Yuen, K.C. (2005). On a correlated aggregate claims model with thinning-dependence structure. Insurance: Mathematics and Economics 36(3), 456-468.

\vspace{0.3 cm}

\noindent Wang, W. and Hu, Y. (2012). Optimal loss-carry-forward taxation for the L\'evy risk model. Insurance: Mathematics and Economics 50(1), 121-130.

\vspace{0.3 cm}

\noindent Wang, W., Wang, Y. and Wu, X. (2018). Dividend and capital injection optimization with transaction cost for spectrally negative L\'evy risk processes. arXiv:1807.11171.

\vspace{0.3 cm}

\noindent Wang, W. and Zhang, Z. (2019). Optimal loss-carry-forward taxation for L\'evy risk processes stopped at general draw-down time. Advances in Applied Probability, 51(3), 865-897.

\vspace{0.3 cm}

\noindent Wang, W. and Zhou, X. (2018). General draw-down based de Finetti optimization for  spectrally negative L\'{e}vy risk processes. Journal of Applied Probability, 55(2), 513-542.

\vspace{0.3 cm}

\noindent Wei, W., Liang, Z. and Yuen, K.C. (2018). Optimal reinsurance in a compound Poisson risk model with dependence.
Journal of Applied Mathematics and Computing 58(2), 389-412.

\vspace{0.3 cm}

\noindent Yao, D., Yang, H. and Wang, R. (2014).
Optimal risk and dividend control problem with fixed costs and salvage value: Variance premium principle.
Economic Modelling 37(1), 53-64. 
\vspace{0.3 cm}

\noindent Yao, D., Yang, H. and Wang, R. (2016). Optimal dividend and reinsurance strategies with financing and liquidation value. ASTIN Bulletin 46(2), 365-399.

\vspace{0.3 cm}

\noindent Yuen, K.C., Liang, Z. and Zhou, M. (2015).  Optimal proportional reinsurance with common shock dependence. Insurance: Mathematics and Economics 64(1), 1-13.

\vspace{0.3 cm}

\noindent Zhang, C. and Liang, Z. (2017). Portfolio optimization for jump-diffusion risky assets with common shock dependence and state dependent risk aversion.  Optimal Control Applications and Methods 38(2), 229-246.

\vspace{0.3 cm}

\noindent Zhao, Y., Chen, P. and Yang, H. (2017).  Optimal periodic dividend and capital injection problem for spectrally positive L\'{e}vy processes. Insurance: Mathematics and Economics 74(1), 135-146.

\vspace{0.3 cm}

\noindent Zhu, J. (2017).  Optimal financing and dividend distribution with transaction costs in the case of restricted dividend rates. ASTIN Bulletin 47(1), 239-268.

\vspace{0.3 cm}


\end{document}